\documentclass{amsart}

\usepackage{amsmath, amsthm}
\usepackage{amssymb, latexsym}
\usepackage{amsfonts}

\usepackage[dvipsnames]{xcolor}

\newtheorem{Theorem}{Theorem}
\newtheorem{theorem}[Theorem]{Theorem}
\newtheorem{Proposition}[Theorem]{Proposition}
\newtheorem{proposition}[Theorem]{Proposition}

\newtheorem{lemma}[Theorem]{Lemma}

\newtheorem{Remark}[Theorem]{Remark}

\theoremstyle{definition}

\newtheorem{alg_step}{}

\newcommand{\thmref}[1]{\rm Theorem~\ref{#1}}

\newcommand{\eqnref}[1]{\rm(\ref{#1})}

\newcommand{\lemref}[1]{Lem\-ma \ref{#1}}
\newcommand{\propref}[1]{\rm Proposition~\ref{#1}}

\newcounter{com}

\newcommand{\beql}[1]{\begin{equation}\label{#1}}
\newcommand{\eeq} {\end{equation}}

    \font\Aaa=msam10

\font\Bbb=msbm10

\newcommand\F{\hbox{\Bbb F}}

\newcommand\B{\hbox{\Bbb B}}

\newcommand\bbb{\hbox{$\mathfrak{b}$}}
\newcommand{\PF}{\mathfrak{P}}

\font\sBbb=msbm8
\newcommand\sB{\hbox{\sBbb B}}
\newcommand\sbbb{{\hbox{\small$\mathfrak{b}$}}}

\font\Aaa=msam10

\def\qed{\hbox{~~\Aaa\char'003}}

\font\Bbb=msbm10

\def\F{\hbox{\Bbb F}}

\DeclareMathOperator{\Tr}{Tr}

        \def\PSL{{\rm PSL}}
        \def\SL{{\rm SL}}

        \def\PGL{{\rm PGL}}

           \def\col{\hspace{-2pt} :  \hspace{-2pt}}

        \def\<{{\langle}}
        \def\>{{\rangle}}

        \font\Aaa=msam10

\font\Aaa=msam10

\font\Bbb=msbm10

\def\F{\hbox{\Bbb F}}

\newcommand{\Fr}{\mathfrak{F}}

\def\a{\alpha}

\def\div{ \kern-.5pt\hbox{\big |} }
\def\ndiv{ {\not\kern-.5pt\hbox{\big |}\,} }
\def\ndivv{ {\not\kern+1.5pt\hbox{$\mid$}\,} }

\def\norm{{\rm N}}
\def\cent{{\rm C}}

\def\diag{{\rm diag}}

\newdimen\refcodesize
\newbox\seriesbox
\def\proof{\noindent {\bf Proof.~}}

\def\PSL  {{\rm PSL }}

\def\SL{{\rm SL }}

\def\SSS{{\mathcal S}}

\def\sbbb{{\hbox{\small$\mathfrak{b}$}}}

\begin{document}

\title[Black box   groups $\PGL(2,2^e)$]
{Black box   groups isomorphic to ${\rm PGL}(2,2^e)$}

\thanks{This research was supported in part by NSA grant MDA-9049810020 and  NSF grants DMS 0900932
and DMS  1303117. }

       \author{W.  M.  Kantor}
       \address{University of Oregon,
       Eugene, OR 97403 and Northeastern University,
Boston, MA 02115}
       \email{kantor@uoregon.edu}
 
    \author{M. Kassabov}
       \address{Cornell University, Ithaca, NY 14853}
       \email{kassabov@math.cornell.edu}


\begin{abstract}
  \vspace{-7pt}
A deterministic polynomial-time   algorithm  constructs an isomorphism between
$\PGL(2,2^e)$ and a black box group
to which it is isomorphic.
\vspace{-9pt}

 \end{abstract}

\maketitle

\centerline{\em Dedicated to the  memory of  \'Akos Seress}

 % \vspace{-7pt}
\section{Introduction}

This   note contains  a polynomial-time algorithm for
\emph{recognizing}
a black box group that is isomorphic to $\PGL(2,2^e)$ by constructing such  an  isomorphism.   
The existence of   a  fast 
algorithm of this type has been open for a number of  years.  The standard way around this existence problem   is based on a lovely idea of Leedham-Green
  \cite{CLG,CLO}, which avoids black box groups, instead
  focussing on groups represented on finite vector spaces   and using  a Discrete Log Oracle
  (this idea applies to all nonzero characteristics).  Unfortunately,
  no black box  version of our result is known in odd characteristic.

 For the required background concerning black box groups,
 including the  timing parameter and straight-line programs (SLPs), see
 \cite{KS,Ser}.  Let
$\mu $   be an upper bound on the time required for
each group operation in  $G=\<\SSS\>$.
We will assume that  $|\SSS|$ is small and hence suppress it  in our timing estimates (not suppressing it would multiply various times by 
 $|\SSS|$).   

The following
 appears to be the first deterministic polynomial-time black box recognition algorithm:

 \begin{theorem}
 \label{Main Theorem}~
 \begin{itemize}
 \item[\rm (i)]
There is an  $O(\mu  e^3\log e)$-time algorithm
which$,$
given a black box group $G=\<\SSS\>$  
isomorphic to $\PGL(2,q)$ with   $q=2^e,$
finds  $3$-element   generating sets $\hat \SSS$ and  $  \SSS^*$ of
$\PGL(2,q)$  and $G,$   respectively$,$
and  a bijection
${\Psi\col \hat \SSS\to    \SSS^*}$
 inducing  an isomorphism
 ${\Psi\col \PGL(2,q)\to G.}$

 \item[\rm (ii)]
In   $O(\mu e^3 )$   time
an {\rm SLP}  of length $O(e)$  
  can be found    from  $\hat \SSS $  to  
    any given   element of  $\PGL(2,q),$  
 and in  $O(\mu e^3 )$   time
an {\rm SLP}  of length $O(e)$  
  can be found    from  $  \SSS ^*$  to 
  any given   element of $G$.
  
   In particular$,$ the  isomorphism 
$\Psi$  is effective$:$ the image of
      any given   element of   $\PGL(2,q)$  and the preimage of    any given   element of $G$ can be found  in   
         $O(\mu e^3 )$  time.

\end{itemize}
\end{theorem}

More precisely,    \emph{$\hat\SSS  $ will essentially be
$\big\{\bigl( \begin{smallmatrix}
1 &0\\
1&1
\end{smallmatrix} \bigr),
\bigl( \begin{smallmatrix}
0&1 \\
1&0
\end{smallmatrix} \bigr), \diag(s ,s^{-1})\big\} \subset \SL(2,2^e),$
where $\F_{2^e}=\F_2[s]$ and the minimal polynomial for $s$ has been computed}
(cf. {\bf \ref{Preliminaries}} and
 {\bf \ref{single}}).  
 (We identify $\PGL(2,2^e)$ and $\SL(2,2^e)$.)
If desired,  one can   switch to  $\F_2[t]$ where  $t$ is a root of 
another 
  irreducible polynomial.

Our new and crucial   observation  is
Proposition~\ref{Krucial}, which can be viewed as providing   
    involution-producing formulas
 in a black box group isomorphic to $\PGL(2,2^e)$
 (cf. Remark~\ref{Las Vegas}).
This is then used to produce  first a Borel subgroup and a
  field $\Fr\cong \F_q$,
and then a group isomorphism.
Also essential is
\propref{calculate entries}, which  recovers the entries of a matrix using the matrices for two given noncommuting involutions.
An  unexpected feature of this result  is that    any given   element of
$G$   can be  obtained quickly using our new generators (\thmref{Main Theorem}(ii)). 
 A  Monte Carlo algorithm related to the Theorem appears in
\cite[Theorem~3.1]{BY}.

Our field calculations all take place inside  $\Fr$, hence essentially
``inside'' $G$,  which explains the timing in
both parts of \thmref{Main Theorem}(ii).
No Discrete Logs are involved, unless a user needs  to
describe the field using a generator of the multiplicative group.
We have not tried to optimize the timing of our algorithm, for example~by using
fast multiplication in $\F_2[s]$ or fast computation of minimal polynomials,
and~we expect that  careful optimization can reduce the time from essentially cubic to  
essentially quadratic in~$e$.%

\section{Preliminaries}
\label{group preliminaries}
The following easy fact is critical:

 \begin{lemma}
 \label{4 points}
If $a,b, c,d$ are distinct points of the projective line over
a field $F,$
 then there is  one and only one  involution in $\PGL(2,F)$ acting as $(a,c)(b,d)\dots$.
 \end{lemma}
 \proof By transitivity we may assume that the points are 
  $\<(1,0)\>, \<(1,1)\>, \<(0,1)\>, $ $\<(s,1)\>, $   
 and then
 $\bigl( \begin{smallmatrix}
0 &1 \\
s&0
\end{smallmatrix} \bigr)$
behaves as required.  \qed
        
 \begin{lemma}
 \label{conjugating involution}
Let $ h\in G=\PGL(2,q)$ with $q>2$ a power of a prime $p$
not dividing  $|h|$.  Let $g\in G$ be such that $[h,h^g]^p \ne1 .$
 Then there is a unique involution  $z\in G$   conjugating
 $h$ to $h^g$.

 \end{lemma}

\proof
If $h$ fixes two points $a,b$ of the  projective line, then  the  hypothesis $[h,h^g]^p \ne1 $ states that
$| \{a,b,c,d\}|=4$ for $\{c,d\}:=\{a,b\}^g$.
By the preceding Lemma
 there is a unique involution $y \in G$ acting as $(a,c)(b,d)\dots$.
Then  $ h^y \in \<h^g\>$ since the stabilizer of 
$c$ and $d$ is cyclic, so that
$  h^y =   (h^g)^\epsilon$  for $\epsilon =\pm1$
since $|\norm_G(\<h\>)\col \<h\>|=2$.

There is also an involution $u$ acting as $(a,b)(c,d)\dots$, and hence commuting with $y$ and satisfying  $  h^{y u} =   (h^g)^{-\epsilon}$.
Then $y$ or $y u$   is the unique involution that conjugates
$h$ to $h^g$.

If $h$ does not fix any points,   embed $G$ into $\PGL(2,q^2)$
and let $\sigma$ be the involutory field automorphism.  Then $h$ fixes two points of the larger projective line.
We have already seen that there is a unique involution $z\in \PGL(2,q^2)$
such that $ h^z=h^g$.  Since $z^\sigma $ also has this property, it follows that $z^\sigma=z $ is in  $G$.  \qed

\begin{Proposition}
\label{Krucial}
Let   $G=\PGL(2,q)$ with $q>2$ even$, $
 let    $2k+1=q^2-1$ denote the odd part of $|G|,$
and let  $1\ne h\in G$  have odd order.
For $g\in G$ with $[h,h^g]\ne 1$ either  $[h,h^g]$  or
 $(hh^g)^kh$    is an involution.

\end{Proposition}
\proof
If  $[h,h^g]$ is not an involution it has odd order.
For $z$ in the preceding  Lemma,  $u:=hz$ has odd order
  (as  $u^2=hh^z=hh^g\ne1$).
Then
 $(h h^g )^kh=(h h^z )^kh=(hz h z )^kh z\cdot z=
(hz )^{2k+1} z=z .      $\qed

\Remark  \rm
\label{Las Vegas}
This result deterministically produces many involutions in a black box group  isomorphic
to $G=\PGL(2,2^e)$.
The unavailability of such a
deterministic or probabilistic
result  has long been the obstacle to a polynomial-time  recognition algorithm
 for $G$.  It has been folklore for several years
 that an involution $u$ would lead to  a Las Vegas algorithm
 based on the following idea:
 find $U:=\cent_G(u) $ using~\cite{Bray,Bo};  use a random element of
 $U$ to produce a random
   element $b$ of $\norm_G(U)$
 (cf.  ``lifting" below  in {\bf\ref{The field}});
 turn $U$ into a field generated by  a single element corresponding to~$b$;~and finally
 make the algorithm Las Vegas
 by verifying  a standard presentation of~$G$.%

\medskip
We will need a standard elementary fact:
\begin{lemma}
\label{conj1}
If $z,z'$ are   involutions in a group and
 $|zz'|$ divides   $2k+1,$
then $(zz')^{k+1}$
conjugates $z$ to $z'$.
\end{lemma}

\section{The square root of a matrix}
\label{The square root of a matrix}
We will use square roots of $2\times2$ matrices over 
$F\cong \F_{2^e},e>1$.
This already occurred in \propref{Krucial}:  $(hh^g) ^{-k}$
 is the square root of
$hh^g$.
  Such square roots    involve very elementary calculations 
   (finiteness is not even required).
We need to exclude the matrices having trace  $0$,
 which   are precisely those of order dividing 2. 

\begin{lemma}
\label{sqrt}
If
$h=\bigl( \begin{smallmatrix}
a&b \\
c&d
\end{smallmatrix} \bigr)$ with $\Tr( h) = a+d\ne0,$ then
$$
\sqrt h =
\frac{1}{\sqrt{a+d}}
 \begin{pmatrix}
a+1&b 
\\
%\displaystyle
c& d+1
\end{pmatrix} 
$$
 is the unique matrix whose square is $h$.
\end{lemma}
\proof
Since $h^2+\Tr(h)h+I=0$ by Cayley-Hamilton, 
$\frac{1}{\sqrt{\Tr(h)}}(h+I)$ is a square root of $h$.  It is unique:
$h$   has odd order $>1$,
 so that  any square root of  $h$ also has odd order and  
hence has to   generate the same
cyclic  subgroup as $h$.
\qed

\medskip
 
For $u=\bigl( \begin{smallmatrix}
1&0 \\
1&1
\end{smallmatrix} \bigr)$,   $ g \in \SL(2, F), |F|=q=2^e$, 
and $q^2-1=2k+1$ as in  \propref{Krucial},  write
 \begin{equation}
\label{B function}
\mbox{$\B(g):= (uu^g)^{k+1}g^{-1}$ ~ if ~
 $[ u,u^g ]  \ne 1 .$}
 \end{equation}
 Since 
$ \sqrt h =h^{k+1}$ in the Lemma, 
 this definition  of the partial function
  $\B$  is based on~\cite{Bray,Bo}
  and can be used for elements of 
  either group in \thmref{Main Theorem}.

%\newpage

\begin{lemma}
\label{calculate B}
If
$g=\bigl( \begin{smallmatrix}
a&b \\
c&d
\end{smallmatrix} \bigr)$ with
$[ u,u^g ]  \ne 1  $  then
 $b\ne0$ and
$\B(g) =
  \begin{pmatrix}
1&0 \\
\displaystyle1+\frac{a+d }{b}&1
\end{pmatrix}  $.
\end{lemma}
\proof
Calculate $h:= {uu^g}  $  and use the preceding Lemma for $\sqrt h$:
$$\hspace{.4in} \B(g)=  \begin{pmatrix}
a&b \\\displaystyle
\frac{ab+a^2+1 }{b}& a+b
\end{pmatrix}
\begin{pmatrix}
d&b \\
c&a
\end{pmatrix}
=
\begin{pmatrix}
1&0 \\ \displaystyle
1+\frac{a+d }{b}&1
\end{pmatrix}. \hspace{.4in}  \qed$$

\medskip
In particular, 
$\B(g)=u$ if and only if $g^2=1$.
Let $\B_{21}(g)$ denote the $2,1$-entry of $\B(g)$.
Starting with knowledge of the entries of $u$, we can use $\B_{21}$ to  find the entries of 
most matrices:

\begin{proposition}
\label{calculate entries}
Let  $r=\bigl( \begin{smallmatrix}
0&\lambda  \\
\lambda ^{-1}&0
\end{smallmatrix} \bigr)$  with $\lambda\in F^*,$
and let $N$ be the normalizer of the maximal torus  containing  $ru$.
If $g\in \SL(2,F)\backslash N,$
then $\lambda$ and the  entries of $g $ can be calculated 
$($by formulas given below$)$ using rational expressions in the square roots of   the
field elements $\B_{21}(g_1gg_2),~g_1,g_2\in\{1,u,r\}^3,$
for which $[ u,u^{g_1gg_2}] \ne 1$.
 
In particular$,$  if  $ \SL(2,F)=\<\SSS\>$ then
$F$ is generated by the set
$$
\big\{
\B_{21}(g_1gg_2) \mid
g\in \SSS \cup \SSS^2, g_1,g_2\in \{1,u,r\}^3,
 \mbox{ for which } [u,u^{g_1gg_2}] \ne 1
\big\}.
$$ 
\end{proposition}
 
\proof
We need to determine the entries of $g=\bigl( \begin{smallmatrix}
a&b \\
c&d
\end{smallmatrix} \bigr)$.
By the preceding Lemma,    \vspace{-4pt}
\begin{equation}
\label{def ABCD}
\begin{array}{r@{\,=\,}l@{\,:=\,}lcr@{\,=\,}l@{\,:=\,}l}
\B_{21}(g)& 1+A& \displaystyle 1+\frac{a+d}{b},
& &\vspace{4pt}
\B_{21}(gr) &1+B & \displaystyle 1+\frac{\lambda ^{-2}b+c}{a},\\
\B_{21}(rg) &1+C & \displaystyle 1+\frac{\lambda ^{-2}b+c}{d}
& \mbox{and} &
\B_{21}(rgr) &1+D & \displaystyle 1+\frac{a+d}{\lambda ^{2}c},
\end{array}
\end{equation}
where $\B$ is not defined in those cases where    the 
denominator is zero. 
Note that $\Tr(g)=\Tr(rgr)=a+d$ and 
$\Tr(rg)=\Tr(gr)=\lambda(\lambda ^{-2}b+c).$
 
If all of the elements $A$, $B$, $C$ and $D$ are   defined and nonzero  (which is equivalent to 
  $abcd(a+d)(\lambda ^{-2}b+c) \not =0$),
then  the relation
\begin{equation}
\label{ABCD}
ABCD\lambda ^2 = (A+D)(B+C)
 \end{equation}
determines $\lambda $.  Moreover,  since $ad+bc=1$,
\begin{equation*}
a=\frac{ACD}{\Delta} , \
b=\frac{D(B+C)}{\Delta}, \
c=\frac{\lambda^{-2}A(B+C)}{\Delta} \ \  {\rm and\ } \
d=\frac{ABD}{\Delta} ,
\end{equation*}
where $\Delta :=\sqrt{\lambda^{-2} AD(B+C)(A+B+C+D) } \not =0$.

If only three of $A,B,C,D$  are defined, then    
$g,gr,rg$ or $ rgr$ has $1,2$-entry  0, so consider the possibility that   $A$ is not defined but $B$, $C$ and $D$ are defined and nonzero (which is equivalent to $acd \not =0$, $b=0$ and $a\not=d$).  Then 
\begin{equation}
\label{BCD}
BCD\lambda ^2 = B+C
\end{equation}
determines $\lambda $.  Moreover,  since $ad=1$,
\begin{equation*}
a=\frac{C}{\Delta} , \
b=0, \
c=\frac{BC}{\Delta} \ \  {\rm and\ } \
d=\frac{B}{\Delta} ,
\end{equation*}
where $\Delta :=\sqrt{BC} \not =0$.
Thus, we can recover $\lambda$ and  the 
 entries of    $g$.  We still must consider the possibility that
  $acd \not =0$, $b=0$ and $a=d$.  Then $ ad+bc=1$ implies that 
 $a=d=1$.  If $c\ne 1, \lambda^{-2}$ then 
 all entries of $urg$ are nonzero   
   and
 $\Tr( ur g) \Tr( r \!\cdot\! urg)\ne 0$, 
 so we can find  $\lambda$ and the entries of 
   $urg$   and hence of $g$ by the previous paragraph.
   If $c = 1$ then $g=u\in N$.  
   If $c=\lambda^{-2}\ne1 $ then 
  $g':= urgu$   satisfies $\Tr(g')  \Tr(rg') \ne0$ and hence has been dealt with already.

If only two of the entries of $g$ are not zero, then $gu$
 has only one zero entry, and we can recover the entries of 
$g$ unless $g\in \{1,r\}\subset N$.  

Finally, we need to deal with the   case where 
$A$, $B$, $C$ and $D$ are all defined and at least one is 0, so that  $abcd\ne0$ and either  $a=d$ or $b=\lambda^2 c$.
Replacing $g$ by $gr$ if necessary, we may 
assume that $a=d$.  
Then   
   $\Tr(gu)= b \ne 0$  and $\Tr(r \!\cdot\! gu )=\lambda(c+a)+\lambda^{-1} b$.
If  $\Tr(r \!\cdot\! gu)\ne0$ then the previous cases determine $\lambda$ and the entries of 
$gu$, hence also of $g$.   On  the other hand, if $\Tr(r \!\cdot\! gu)=0$ then 
$a=c+\lambda^{-2 }b$, which is precisely the condition that the 
involution $g$ inverts $ru$ and hence lies in $N$.
 (Note that it is not surprising that it is 
not possible to recover  $\lambda$ if $g\in N\backslash \{1,u\}$,
since  then   $\B(g) \in N \cap \cent_G(u) = \{ 1, u \}$ by   \eqref{B function}.)

For the last statement of the proposition, note that  the entries of 
 $\SSS$   generate $F$. 
We have already computed the entries of every 
generator $s' \in \SSS  \backslash  N$,
while
 we can  obtain the entries of any  $s\in \SSS \cap N$  from the entries of   $s'$ and $ss'$.
(In fact we do not need  $\SSS^2$, only   
$\SSS\cup s'\SSS$ 
 for a single 
$s'\in \SSS\backslash N$.)
\qed

\Remark  \rm The  calculation of $\B(g)$
in  \lemref{calculate B} was unexpectedly
simple.  Although it will be less simple in larger groups having more elements of even order,
such a calculation might occasionally  be useful  in order to speed up
the recovery of the entries of a matrix from various values of $\B$
(as in \propref{calculate entries} and
{\bf\ref{find g}} below).

%\vspace{-1pt}
\section{Proof of Theorem~\ref{Main Theorem}}%({\rm i})}
We are
given a black box group $G=\<\SSS\> \cong \PGL(2,q)$.
We may assume that $q=2^e>2$. 
We proceed in several steps, each of which 
begins with a hint of its content.

\begin{alg_step}\rm
\label{Step 1}  (Noncommuting involutions $u,r$.)
If every element of $\SSS$ is an involution then two of them do not commute and we  obtain noncommuting involutions $u,r$.
If some $1\ne h\in \SSS$ has odd order,
 then 
 $ \cent_G(h) $ is cyclic and not normal in $G$,  and hence  $h^g\notin  \cent_G(h)$
 for some $g\in \SSS$, so that
 Proposition~\ref{Krucial}
produces an involution $u$.
Then 
some $s\in \SSS$ does not normalize $\cent_G(u) $, so that 
 $r:= u^s\notin  \cent_G(u) $.
(Time:  $O(\mu e)$ using SLPs of length $O (e)$.)

\end{alg_step}

%\newpage

\begin{alg_step}\rm
\label{The field: motivation}
(The field: motivation.)
In order to define a field we need to understand consequences of the assumed isomorphism $\hat G:=\PGL(2,F)\to G,$ 
where $F\cong \F_q$. 
(In  {\bf \ref{The field}} 
we will define  a specific model of $\F_q$.)

Let 
$\hat u:=\bigl( \begin{smallmatrix} 1 & 0\\ 1 &1 \end{smallmatrix} \bigr) $, $\hat U:=\cent_{\hat G}(\hat u)$,  
$\hat B:=\norm_{\hat G}(\hat U)$  and $\hat B^*:= \hat B/\hat U$.
Writing matrices modulo 
scalar matrices,
 $\hat B^* $ consists of cosets 
 ${\tilde s}\hat U=\bigl( \begin{smallmatrix} s  & 0\\ * & 1 \end{smallmatrix} \bigr) $ 
 that can be viewed as field elements.
For such cosets ${\tilde s}\hat U $ and 
${\tilde t}\hat U , $ 
their product in $F^*$ 
corresponds to the group operation on 
$\hat B^*$, while   
 field-addition   
occurs via\vspace{2pt}
 $\hat u^{{\tilde s}+{\tilde t}} =\hat u ^{\tilde s}\hat  u^{\tilde t}$
 when  ${\tilde s}\hat U\ne {\tilde t}\hat U$.
 
 Each   $\hat x \in \hat  U$ has the form $\hat x =\hat X(t):=  
 \bigl( \begin{smallmatrix} 1 & 0\\ t &1 \end{smallmatrix} \bigr)$,
 $t\in F$.  
 \emph{If  $\hat r
 =\bigl( \begin{smallmatrix}
0&\lambda  \\
\lambda ^{-1}&0
\end{smallmatrix} \bigr)$  with $\lambda\in F^*$ 
\emph{$($as in 
{\propref{calculate entries}}$) $} and $t\ne0, $ then}
(using     $2k+1=q^2-1$)
  \begin{equation}
\label{hat b}
\mbox{
$\bbb(\hat x):=
 (\hat u\hat u^{\hat r} )^{k+1} ( \hat u^{\hat r}\hat x )^{k+1} \hat U
 = \bigl( \begin{smallmatrix}
 \sqrt t & 0 \\ 0& \sqrt t^{\,-1} 
  \end{smallmatrix} \bigr)  \hat U
  \in \hat B/\hat U $,
 }
 \end{equation}
 which is independent of $r$.
  For,   $\hat u\hat u^{\hat r } $ and $\hat u^{\hat r}\hat x$    have odd order (since $\hat u^{\hat r}\notin\hat U$).  Then
 $ %\bbb(\hat x):=
 (\hat u\hat u^{\hat r}  )^{k+1} ( \hat u^{\hat r}\hat x )^{k+1} $  conjugates $\hat u$ to $\hat x$ (by   \lemref{conj1}), 
  hence normalizes $\hat U$
and so lies in $\hat B$.
A simple matrix calculation in $\hat B$ produces all matrices conjugating $\hat u$ to $\hat x$,  as stated
  in
\eqref{hat b}.

\end{alg_step}

\begin{alg_step}
\label{The field}
(The field $\Fr$.)
We   define what amounts to a ``black box field''
$\Fr$ obtained from  the involution $u$.  Define 
$ U:=\cent_{  G}(  u)$,  
$  B:=\norm_{  G}(  U)$  and $  B^*:=   B/  U$.
Then $\Fr$ arises from the pair consisting of
$U $ and $B^* $ as additive and multiplicative groups, respectively;
as a set $\Fr$ will be $B^*  \cup \{0\}  $ with $0$ treated in the obvious ways.  We do not yet  have   these as constructed groups.
%\marginpar{constructed groups}

\begin{itemize}
\item
If $g\in G$,  then 
$g\in U\iff [u,g]=1$, and $g\in B\iff [u,u^g]=1$.

\item If  $t_1U,t_2U\in B^*$ then  $t_1U=t_2U \iff(t_1t_2^{-1})^2=1$.

\item  \emph{Lifting elements of $U$ to $B/U$}:
If $  1\ne x\in U $, then
 $uu^r  $ and $u^rx$    have odd order (since $u^r\notin U$),  so that
(as in \eqref{hat b})
\begin{equation}
\label{b}
\mbox{
$\bbb(x):=
 (uu^r  )^{k+1} ( u^rx )^{k+1} U$ \  conjugates $u$ to $x$.
 }
 \end{equation}
Then
 $\bbb(x)$    normalizes $U$ and  so lies in $B/U$.
(Time: $O(\mu e)$  using SLPs of length $O (e)$.)

\item \emph{Labeling $ U\backslash\{1\}$ using $\Fr^*=B^*$}: If $tU\in B^*$ then
$X(tU):=u^t \in  U\backslash\{1\}$
defines   the inverse of the map
$\bbb\col U\backslash\{1\} \to B^*$
in  \eqnref{b}.  Also let  $X(0) :=1$.

\item \emph{Field multiplication}:~This  is inherited from $B/U$ and hence from $G$.
(Time:~$\mu$.)

\item  \emph{Field addition}:
$t_1U+t_2U = tU=\bbb(u^{t_1} u^{t_2})$
for distinct $t_1U,t_2U,$ so that
 $u^t=u^{t_1} u^{t_2}$.~(Time:~$O(\mu e)$.)
 \end{itemize}

All of the above imitated {\bf \ref{The field: motivation}}.  
We emphasize that the field depends entirely on $u$ 
(which uniquely determines $U$ and $B$):
it does not depend on   $r$, which was used only to obtain 
 elements of $B^*=B/U$.
 {We also emphasize that nonzero field elements are cosets:
whenever we write a nonzero field element we are implicitly also writing (and storing) a coset representative.}

Field-theoretic calculations  are   postponed to Section~\ref{Field computations}.    \
\end{alg_step}

\begin{alg_step}\rm
\label{generate field}
(Generators of $\Fr$.)
\emph{The elements
$\bbb(\B(g_1g g_2))$ with 
 $ g \in \SSS \cup \SSS^2 ,g_1,g_2\in \{1,u,r \}^3$
and $[u,u^{g_1 g g_2}]\ne1,$
generate
  $\Fr \cong \F_q$.}
Here the partial function  $\B$ is defined as
 in \eqnref{B function},  $\bbb$ was defined in
\eqref{b},
  and  ``generate''   means that the
stated elements  lie in no proper subfield.
(Time:  $O(\mu  e)$ to find these elements
using SLPs of length $O(e)$.)

For,  there is  an isomorphism  $\Psi\col  \PGL(2,\Fr)\to G$.
By a basis  change  of $\Fr^2$ we may assume that  $\Psi$
sends $\hat u:=
\bigl( \begin{smallmatrix}
1&0  \\
1&1
\end{smallmatrix} \bigr)
\mapsto u$ and
$\hat r:=
\bigl( \begin{smallmatrix}
0&\lambda^{-1}  \\
\lambda&0
\end{smallmatrix} \bigr)  \mapsto r$  for some uniquely determined
$\lambda\in \Fr$.  
Since the preimages of
$\SSS$ cannot all lie in a proper subfield of $\Fr $,
  the   elements 
$\bbb(\B(g_1 g g_1) )$ arising  from  the last part of 
  \propref{calculate entries}
do not lie in a proper subfield
of  $\Fr $.

 %\newpage

 \Remark
 \label{isomorphism}
  \rm The isomorphism $\Psi$ sends $\hat u\mapsto u$ and  induces an isomorphism
from the usual Borel subgroup $\hat B$  of $\hat G$ to $B$.  The latter isomorphism is uniquely determined up to conjugation and field automorphisms.  Hence,
we may assume in our arguments that
 $\hat X(t)  \mapsto X(t)$.
(We emphasize \emph{in our arguments} since we have yet to construct such an isomorphism, as required in the theorem.)

We will abbreviate 
$\Psi':=\Psi^{-1}$  and $\hat g := \Psi'(g)$.

%\color {blue} 

\Remark  \rm (The projective line.)
\label{projective line}
  We also obtain the projective line $\PF = \Fr \cup \{\infty\}$.
  The fact that we   can obtain   the  action of  any $g\in G$ on it
  motivated parts of this paper.
It  suffices  
to show how to compute $0^g$ and $\infty^g$.
For, any $x \in \Fr^*$  is represented by an element $b\in B$,
so that  $x^{b^{-1}u} =0$ by definition and hence   $x^g = 0^{ubg}$.

By the proof of Proposition~\ref{calculate entries},
if 
$\B(g)$, $\B(gr)$, $\B(rg)$ and $\B(rgr)$ are all defined and
are not equal to $u$ (so that $g\notin N$), then
$$
\begin{array}{r@{\,}c@{\,}l}
0^g &=& \bbb\big(u\B(g)\big) \bbb\big(u\B(gr)\big) \Big(\bbb\big(\B(g)\B(rgr)\big) \Big)^{-1} % \quad \mbox{and}
 \vspace{.01pt}
\\
\infty^g &=& \bbb\big(u\B(g)\big) \bbb\big(u\B(gr)\big) \Big(\bbb\big(\B(rg)\B(gr)\big) \Big)^{-1}.
\end{array}
$$

\end{alg_step}

\begin{alg_step}\rm
  \label{r'}
  
  ($\PGL(2,2)$ and  $r'$.)
  By the proof of \propref{calculate entries},    there 
is an element 
%\color{red}
 $g\in  \{1, r, u  \}^3  \SSS  \{1, r, u  \}^3$
   such that three or four of the  elements
  $\B(g)$,  $ \B(gr),\B(rg), \B(rgr)$ are defined,  and each of those is not equal to $u$.
   
 {\em   If all  four  behave this way$,$
calculate
$$
\tau : =
\bbb\big(\B(g) \B(rgr)\big)  
\bbb\big(\B(rg)\B(gr)\big) 
\Big(
\bbb\big(\B(g)u\big) 
\bbb\big(\B(rg)u\big)
\bbb\big(\B(gr)u\big)
\bbb\big(\B(rgr)u\big)
\Big)^{-1}  \! \!\!.
$$
Let $\tau=tU,  t\in B$.  
Then $uu^{tr u}$ has order $3$ and
$r':=u^{tru}$ is an involution$,$
so that $\<u,r' \>\cong \PGL(2,2)$.}  (Time:  $O(\mu e)$.)

For,    the definitions of $\B$ and $\bbb$
(in   \eqref{B function},   
\eqref{hat b} and \eqref{b})
  imply that  
  $$\hat u^{\Psi' ( \sbbb(\sB(g) \sB(rgr)  ) ) }=
 \Psi'(u^{ \sbbb(\sB(g)\sB(rgr)) })=
 \Psi'(\B(g)\B(rgr) ) =
  \B(\hat g)\B(\hat r\hat g\hat r) , $$
 and hence  that
   $\Psi'\big( \bbb\big(\B(g)\B(rgr) \big) \big)= 
  \bbb\big(\B(\hat g )\B(\hat r\hat g\hat r) \big)$;
 and similarly for the other terms in the definition of $\tau$.
  Using  \eqref{def ABCD}, \eqnref{ABCD}
  and \eqref{hat b},
we obtain
$$
\begin{array}{lllll}
\Psi ' (tU)  \hspace{-6pt}  &= \hspace{-6pt}&
\bbb\bigl( \begin{smallmatrix} 1&0  \\  A+D &1 \end{smallmatrix} \bigr)
\bbb\bigl( \begin{smallmatrix} 1&0  \\  B+C &1 \end{smallmatrix} \bigr)
\Big(
\bbb\bigl( \begin{smallmatrix} 1&0  \\  A &1 \end{smallmatrix} \bigr)
\bbb\bigl( \begin{smallmatrix} 1&0  \\  B &1 \end{smallmatrix} \bigr)
\bbb\bigl( \begin{smallmatrix} 1&0  \\  C &1 \end{smallmatrix} \bigr)
\bbb\bigl( \begin{smallmatrix} 1&0  \\  D &1 \end{smallmatrix} \big)
\Big)^{-1}\vspace{4pt}
\\
%\bbb
&= \hspace{-6pt}&
\left( \begin{smallmatrix} 
\sqrt{\textstyle\frac{(A+D)(B+C)\rule[-0.1ex]{0pt}{1.7ex}}{ABCD\rule[1.7ex]{0pt}{0pt}}} & 
\textstyle 0  \\  
\textstyle 0 &  
\sqrt{\textstyle\frac{(A+D)(B+C)\rule[-0.1ex]{0pt}{1.7ex}}{ABCD\rule[1.7ex]{0pt}{0pt}}}^{\scriptstyle\,\, -1\!\!}
\end{smallmatrix} \right)
 \hat{U} \, = \,
\Big( \begin{smallmatrix} \lambda &0  \\  0 & \lambda^{-1} \end{smallmatrix} \Big) \hat{U}.
\end{array}
$$
Since $\hat r=
\bigl( \begin{smallmatrix}
0&\lambda^{-1}  \\
\lambda&0
\end{smallmatrix} \bigr)  $, it  follows that 
$\widehat{r'}=\Psi'(r')=\hat u^{\Psi'(tU)\hat r \hat u}
= 
{\hat u \rule[2.5ex]{0pt}{0pt}}^{
\left(\! 
\begin{smallmatrix} \scriptscriptstyle 0 & \scriptscriptstyle 1  \\  \scriptscriptstyle 1 & \scriptscriptstyle 0 \end{smallmatrix}
\!\right)
\left(\! 
\begin{smallmatrix} \scriptscriptstyle 1 & \scriptscriptstyle 0  \\  \scriptscriptstyle 1 & \scriptscriptstyle 1 \end{smallmatrix}
\!\right)
}
=\bigl( \begin{smallmatrix} 0 &1  \\  1 &0 \end{smallmatrix} \bigr)
$,
so that $\<u,r'\>\cong \<\hat u,   \widehat{r'}\>\cong \PGL(2,2)$.
 
Similarly, assume that  exactly
three   of the  elements
%{\color{red}  
$\B(g)$, $\B(gr)$, $\B(rg)$, $\B(rgr)$ are defined  and not equal to $u$.
\emph{If   $  \B(g)$ is not defined$,$  then 
  {\eqnref{BCD}}  can be used   as above in place of 
    {\eqnref{ABCD}$,$}
    this time together with  the field element      
    $$\tau   : =   
 \bbb\big(\B(rg)\B(gr)\big) 
\Big(
\bbb\big(\B(rg)u\big)
\bbb\big(\B(gr)u\big)
\bbb\big(\B(rgr)u\big)
\Big)^{-1}  \! \!\! , $$  
letting $\tau=tU$  and
$r'=u^{tru}$
in order to obtain 
$\<u,r' \>\cong \PGL(2,2)$.}
The remaining cases are handled by replacing $g$ by $rg$, $gr$ or $rgr$. 
  
\medskip
 Now we can introduce
  some standard elements of $\hat G$ and $G$  for  all  $t\in \Fr^*$:
        $$
 \mbox{$\hat n(t) := \hat X(t)\hat X(t^{-1})^{\widehat{r'}}\hat X(t)$  and   $\hat h(t) := \hat n(t)\hat n(1)$}
 $$
  $$
 \mbox{ $n(t) := X(t)X(t^{-1})^{r'}X(t)$  and   $h(t) := n(t)n(1)$.}
 $$
 Then $\Psi $ sends  $ \widehat {r'}
 \to r'$ and
 $\hat h(t)   \mapsto h(t). $
\end{alg_step}

  \begin{alg_step}\rm
{\bf Proof of \thmref{Main Theorem}(i).}
In  {\bf \ref{single}} we will express $\Fr$ in the form $\F_2[s]$.
Since the  isomorphism $\Psi$  sends
$\hat X(1) \mapsto  X(1),   \widehat {r'}  \mapsto{r'}
$ and    $\hat h(s) \mapsto h(s)$,
  the map
 $  \hat \SSS:=\{ \hat X(1), \widehat {r'} , \hat h(s) \}
  \to \SSS^*:=\{  X(1),    {r'},  h(s) \}$
determines $\Psi$.

The stated time  is dominated  by {\bf \ref{single}}.
\end{alg_step}

\begin{alg_step}
\label{find g}
{\bf Proof of \thmref{Main Theorem}(ii).}
In order to handle elements of $\hat G$, 
first consider  $\hat X(t)$, $t\in \Fr$.  In $O(\mu e^3 )$ time
write $t$ as  $\sum_0^{e-1}a_is^i$ with $a_i\in \F_2$
  (cf.~{\bf \ref{Linear algebra}}),
and then find  an SLP
from  $\hat \SSS$  
 to $\hat X(t)=\prod_i\hat X(s^i)^{a_i}$
of length $O(e)$.
Now it is
 easy to use the definitions at the end of 
 {\bf \ref{r'}}, together with 
straightforward row reduction,
      to  find   an SLP  of length $O(e)$ from $\hat \SSS$  to  any
      given  element in $\hat G = \hat B\cup \hat B\hat r\hat B$, 
      where  $\hat B=\hat U\{\hat h(t)\mid t\in \Fr^*\}$.
 This takes $O(\mu e^3)$ time
 (dominated by~{\bf \ref{Linear algebra}}).
  
Now consider $g\in G$.
We must find an SLP from $\SSS^*$ to $g$.
 Assume that $g\notin \norm_{  G}(\<  u ,   {r}\> ) $.
For each $ g_1, g_2\in  \{1,u , r\}^3  $
such that $[u,u^{g_1 g g_2}]\ne1$,  find
$\B(g_1 g g_2)\in   U$ and 
$\B(\hat  g_1\hat   g\hat   g_2)
 = \Psi'(\B(g_1 g g_2))\in \hat  U$, and then use
 \propref{calculate entries} to find  the entries of the matrix $\hat  g = \Psi'(g)$.  
 Use the preceding paragraph to find an SLP from 
   $\hat \SSS$ to  $\hat g$  of length $O(e)$.
     Then the $\Psi$-image of that SLP  is the required SLP
     from $\SSS^*$  to $g$.
Finding $\hat g$ takes  $O(\mu e)$ time,  finding an SLP 
to  $\hat g$ takes
  $O(\mu e^3  )$ time,
  and  applying $\Psi$ to the members of the SLP  takes
  $O(e  )$ time.

 If $g\in \norm_{  G}(\<  u ,   {r}\> )  $ then
 $g=g h(s)^{-1}\! \cdot \! h(s)$ with $h(s)\in \SSS^*$,
 and we can find an SLP from $\SSS^*$ to $g h(s)^{-1}\notin 
 \norm_{  G}(\<  u ,   {r}\> ) $ in the required time.
  \qed
  
\end{alg_step}

%\newpage

\section{Field computations  }
\label{Field computations}

 The proof of  \thmref{Main Theorem} depends on
   the black box field $\Fr$  in
  {\bf \ref{The field}} that   is a model of the field $\F_q$ inside the black box group $G$. Using it we will  emulate 
some  classical field-theoretic algorithms.
 Recall    that each addition in $\Fr$ takes $O(\mu e)$ time while  each multiplication takes
$O(\mu  )$ time.

%\medskip
\begin{alg_step}  (Preliminaries.)
\label{Preliminaries}
{\rm(i)}
Given $s\in \Fr^*,$   find $m$ such that $ \F_2[s]\cong \F_{2^m}$
in $O(\mu  e  )$ time: find the first $m\in \{1,\dots,e\}$
such that $s^{2^m}=s$.
In particular, we can test whether  $ \F_2[s]= \Fr .$

{\rm(ii)} For $s$ and $m$ in (i), the minimal polynomial of $s$ over $\F_2$ is $\prod_0^{m-1}(x-s^{2^i})$, found in   $O(\mu e^3 )$  time.
This yields an isomorphism    $\F_2[x]/(f)\to \F_2[s] $   induced by
$x+(f)\mapsto s$.

{\rm(iii)}  The trace map $\Tr\col \Fr\to \F_2$, defined by
$\Tr(t):=\sum_0^{e-1}t^{2^i}$, is  calculated  in
$O(\mu   e^2  )$ time
using
$u^{\Tr(t)}=\prod_0^{e-1} u^{t^{2^i}  }$.

{\rm(iv)}   In
 {\bf \ref{Linear algebra}} we will use  linear equations over $\F_2$.  While this $\F_2$ is contained in $\Fr$, with its rather slow field operations,
we can work with a standard model of $\F_2$ with much faster field operations.

\end{alg_step}

%\medskip
\begin{alg_step}
(Linear algebra.)
\label{Linear algebra}
Let  $s\in \Fr$ with $\Fr=\F_2[s]$.
\emph{In $O(\mu  e^3 )$ time$,$ given $t\in \Fr$
we  can  find 
the unique $x_i\in \F_2$ such that
\vspace{1pt}
$t=\sum_0^{e-1}x_is^i$.}

First find all  $2e-1 $ traces   $a_{ij}:=\Tr(s^{i+j})  \in \F_2$
and  all
$\Tr(ts^j) \in \F_2,$  $0 \le i,j <e$, in
$O(\mu  e^ 3)$ time.
Since
$\Tr(ts^j)=\sum_{i=0}^{e-1}x_i\Tr(s^{i+j})$,
we obtain
$$\vspace{-2pt}
\Tr(ts^j)=\sum_{i=0}^{e-1}a_{ij}x_i, ~~0\le j<  e.
\vspace{-1pt}
$$
In $O( e^3)$ time 
 solve these linear equations over $\F_2$.      %

     \end{alg_step}

%\newpage

\begin{alg_step}
\label{single}
(Field generator.)
\emph{In $O(\mu e^3 \log e)$ time we can find
$s$ such that  $\Fr=\F_2[s]$.}

By   {\bf \ref{generate field}},  $ \F_{2^e}\cong  \Fr=\F_2[\a_1,\dots,\a_n]$ with  $n=O(1)$
(we are ignoring $ |\SSS|$).

1. Factor  $e = \prod _i q_i$  into  powers $q_i=p_i^{k_i}$
of  $O(\log e)$ different primes $p_i$.  (Time:  $O(e)$.)

2. Find the degree of each $\alpha_j$  over $\F_2$.
(Time:  $O(\mu e^2)$ as in
 {\bf \ref{Preliminaries}}(i).)

3. For each $i$  find $ \alpha_{j_i}$  such that $q_i $ divides
$ m_i:= [ \F_2[\alpha_{j_i}] \col  \F_2]$ (this exists since the
$ \alpha_j $  generate   $\Fr$).

4. For each $i$ compute the  polynomial
$f_i  : = \prod_{t=0}^{(m_i /q_i)-1}(x-\a_{j_i}^{2^{t q_i}})$,
whose coefficients generate a field  $ \F_{2^{q_i}}. $
 (Time:  $O(\mu e^3 \log e )$.)

5.  For each $i$   find a coefficient $c_i$ of $f_i$ such that  $\F_{2^{q_i}}=\F_2[c_i] $      (this exists since $q_i$ is a prime power).
 (Time:  $O(\mu e^2\log e)$, testing all coefficients using
   {\bf \ref{Preliminaries}}(i).)

6.  Output $s:= \prod _i c_i$.  Then $ \F_q=\F_2[s] $
 since the groups $ \F_2[{c_i}]^* $ have pairwise  relatively prime orders, so that
 each $c_ i$ is a power of $s$.
(Time: $O(\mu  \log e)$.)

\end{alg_step}

At this point we have completed the requirements made after the statement of
\thmref{Main Theorem}:  we have obtained $\Fr$ as $\F_2[s]$, and we
can find the minimal polynomial of $s$ (cf.  {\bf \ref{Preliminaries}}(ii)).

\vspace{-2pt}
\section{Odd characteristic}
We conclude with  remarks about a case we have not been able to  handle:
$G=\PSL(2,q)$ with $q=p^e$ for an odd prime $p$.
Since about a quarter  of all group elements have even order, it is easy to
(probabilistically)  find an involution.  The problem is to find an element of order $p$.

There is an analogue of Proposition~\ref{Krucial} that might not be entirely useless, though we do not see how to use it.
That result can   be imitated given an element $h$ such that $|hh^g|$ is  a factor of the odd integer
$(q\pm1)/2$;  if the involution $t$ in Lemma~\ref{conjugating involution} is in $\PSL(2,q)$, then we  obtain that involution as before.

There is an analogue of Lemma~\ref{4 points} involving elements of order $p$ acting as $g=(a,c,\dots)(b,d ,\dots)\dots$:
given distinct $a,b$, for about half of all pairs $c,d$ there is such an element
$g$ (in fact, two of them), and then $G_{ab}^g=G_{cd}$.
 However, we do not see how to use this fact.

\medskip
\noindent
\emph{Acknowledgement}:
We are grateful to
Peter Brooksbank for helpful suggestions and for implementing 
our algorithm in Magma.
We thank the referee
 for his statement of \lemref{sqrt}.

%\vspace{-4pt}

\end{document}